\def\no{\if01}
\def\iftwelvept{\no}

\def\ifusepdf{\no}
\def\ifpsfont{\no}

\iftwelvept
\documentclass[leqno,12pt]{amsart}
\else
\documentclass[leqno,12pt]{amsart}
\fi
\usepackage{amssymb}
\usepackage{amscd}
\usepackage{latexsym}
\usepackage{verbatim}
\usepackage{mathrsfs}
\usepackage[all]{xy}

\ifusepdf
\usepackage{hyperref}
\else\fi
\ifpsfont
\usepackage[T1]{fontenc}
\usepackage{times}
\else\fi

\iftwelvept
\else\fi

\theoremstyle{plain}
\newtheorem{Theorem}{Theorem}[section]

\newtheorem{Proposition}[Theorem]{Proposition}

\theoremstyle{definition}

\newtheorem{Remark}[Theorem]{Remark}

\renewcommand{\theTheorem}{\arabic{section}.\arabic{Theorem}}

\newcommand{\ZZ}{{\mathbb{Z}}}
\newcommand{\QQ}{{\mathbb{Q}}}

\newcommand{\YY}{{\mathcal{Y}}}

\newcommand{\OO}{{\mathcal{O}}}

\newcommand{\Spec}{\operatorname{Spec}}

\newcommand{\XX}{\mathcal{X}}

\newcommand{\Proof}{{\sl Proof.}\quad}
\newcommand{\QED}{{\unskip\nobreak\hfil\penalty50\quad\null\nobreak\hfil
{$\Box$}\parfillskip0pt\finalhyphendemerits0\par\medskip}}

\begin{document}
\title{Note on local structure of Artin stacks}
\thanks{The author is supported by JSPS}
\author{Isamu Iwanari}

\address{Research Institute of the Mathematical Sciences,
Kyoto University, Kyoto, 606-8502, Japan}
\email{iwanari@kurims.kyoto-u.ac.jp }

\maketitle

\renewcommand{\theTheorem}{\arabic{Theorem}}

Let $\XX$ be an Artin stack over $\ZZ$ (see \cite{LM}).
A coarse moduli map (space) for $\XX$ is
a morphism $\pi:\XX\to X$ to an algebraic space $X$ which has
properties:
\begin{itemize}
\item $\pi$ is universal among morphisms from $\XX$ to algebraic spaces,
i.e., for any morphism $\phi:\XX\to Y$ with an algebraic space $Y$ there is
a unique morphism $f:X\to Y$ such that $\phi=f\circ \pi$,
\item for any algebraically closed field $K$, $\pi$
identifies the set of isomorphism classes of $\XX(K)$ with
the set $X(K)$ of $K$-valued points of $X$.
\end{itemize}

Let $\XX$ be an Artin stack locally of finite type over
a locally noetherian scheme $S$.
Suppose that $\XX$ has finite inertia stack, that is, the projection
\[
\textup{pr}_1:I_{\XX}:=\XX\times_{\XX\times_S\XX}\XX\longrightarrow \XX
\]
is finite.
In this setting, according to \cite{KM}
there exist an algebraic space $X$ locally of finigte type over $S$
and a coarse moduli map $\pi:\XX\to X$ over $S$, where
$\pi$ is proper and quasi-finite.
When in addition $\XX$ is a Deligne-Mumford stack,
it has, \'etale locally on its coarse moduli space $X$,
the form of a quotient stack $[Z/G]$ where $Z$ is an affine scheme
and $G$ is a finite group. This structure has been important
and very useful in various situations.
We would like to make the following useful observation:

\begin{Theorem}
\label{main}
Let $\XX$ be an Artin stack locally of finite type over $S$ with finite
inertia stack.
Then for any point $x\in \XX$ there is an \'etale neighborhood
$U\to X$ of $\pi(x)$ such that the $S$-stack $\XX\times_XU$ has 
the form $[\Spec R/\operatorname{GL}_n]$, where
$\Spec R$ is an affine scheme over $S$ and $\operatorname{GL}_n$
is the general linear group scheme defined over $S$.
Moreover, $U$ is naturally isomorphic to the
spectrum $\Spec R^{\operatorname{GL}_n}$ of the invariant ring.
\end{Theorem}

\Proof
Let $\pi(x)=y$.
Then according to the proof of the theorem of Keel and Mori (see \cite{KM}
in particular Section 4 in loc. cit.), we see that
there is an \'etale neighborhood
$U\to X$ of $y$ from an affine scheme $U$ such that
$\XX_U:=\XX\times_XU$ has a finite flat surjective morphism $p:Y\to \XX_U$
from a scheme $Y$.
Then since $Y$ is a scheme, $p_*\OO_Y$ is a vector bundle on $\XX_U$
such that at each geometric point
the stabilizer faithfully acts on the fiber.
After shrinking $\XX_U$ if necessary we assume that $\XX_U$ is connected.
Then by \cite[2.12]{EHKV}
the total space of
the corresponding principal $\operatorname{GL}_n$-bundle
$W$ over $\XX_U$ is an algebraic space
and thus $\XX_U$ is isomorphic to the quotient stack
$[W/\operatorname{GL}_n]$ where $n$ is a non-negative integer.
Now we will show that $W$ is an affine scheme.
Note that the projection $W\to [W/\operatorname{GL}_n]$
is a $\operatorname{GL}_n$-bundle and $\XX_U\simeq[W/\operatorname{GL}_n]\to U$
is proper. Hence $W$ is separated and noetherian.
In addition, $Y\to \XX_U\to U$ is a finite morphism and thus
$Y$ is affine.
Therefore $W\times_{\XX_U}Y$ is affine because the
second projection $W\times_{\XX_U}Y\to Y$ is affine.
Notice that the first projection $W\times_{\XX_U}Y\to W$
is a finite surjective morphism.
Consequently, applying Chevalley's theorem for algebraic
spaces (see \cite[III 4.1]{Kn})
we deduce that $W$ is affine.
The last claim follows from the fact that $\XX_U\to U$ is
a coarse moduli map because $U\to X$ is \'etale.
\QED

\begin{Remark}
Suppose further that $\XX$ is {\it normal} in Theorem~\ref{main}.
Then we can take an action of $\operatorname{GL}_n$ on $\Spec R$
in Theorem~\ref{main}
to be linearized.
Note that under the assumption $\Spec R$ is normal.
According to \cite[Theorem 3.3, 2.18]{Th} (see also \cite[2.5]{Su}),
there is a $\operatorname{GL}_n$-vector bundle $\mathcal{V}$ on $U$ and a $\operatorname{GL}_n$-equivariant
immersion $\Spec R\to \mathbb{P}(\mathcal{V})$.
\end{Remark}

\begin{Remark}
If $\XX$ is a Deligne-Mumford stack of finite type over $S$
which has finite inertia stack, then \'etale locally on its coarse moduli
space, $\XX$ is the quotient $[Y/G]$ of an affine scheme $Y$
by an action of a finite
(constant) group $G$ (see \cite{AV}).
If stabilizer group schemes at geometric points on $\XX$
are finite (not necessarily reduced) linearly reductive group schemes,
in \cite{AOV} it is shown that $\XX$ is, \'etale locally on its coarse moduli space, the quotient
of an affine scheme by an action of a linearly reductive group scheme.
In \cite{Is} such \'etale-local quotient structures
was studied when $\XX$ has (not necessarily finite) linearly reductive
stabilizers and satisfies the stability (see \cite{Is}).
In these cases the stabilizer group schemes
has
no non-trivial deformation (i.e., has a unique deformation),
and we may take $G$ to be the stabilizer
group scheme at a point on $\XX$
when we work over a field.
This point is crucial for Luna's \'etale slice theorem.
However, in positive characteristic case a general (finite) group scheme
has many and rich deformations. Thus in general local structures of
Artin stacks are not so simple as above cases. Indeed
an Artin stack can contain the information arising from
non-trivial flat deformations of $G$, i.e., $BG$.
\end{Remark}

The typical usage of Theorem~\ref{main}
is the reduction of problems to the case of group actions.
In the rest of this note, we will present one of such applications
of Theorem~\ref{main}, which is a direct one.
For this we shall prepare our setup.

Let $\XX$ be an Artin stack of finite type over $\ZZ$.
We will denote by $G(\XX)$ (resp. $K(\XX)$)
the algebraic K-theory spectrum of the exact category of
coherent sheaves (resp. vector bundles) on $\XX$, and we
let $G(\XX)\otimes\QQ$ and $K(\XX)\otimes\QQ$
Bousfield localizations of $G(\XX)$ and $K(\XX)$ respectively,
with respect to $\QQ$ (see \cite{Jo}).
Let us recall the isovariant \'etale descent of $G$-theory
due to Joshua \cite{Jo}, which generalizes Thomason's descent \cite{Thae}.
A morphism $\YY\to \XX$ of Artin stacks is isovariant
if $I_{\YY}\to I_{\XX}\times_{\XX}\YY$ is an isomorphism,
where $I_{\XX}$ and $I_{\YY}$ denote inertia stacks of $\XX$ and $\YY$
respectively. Note that an isovariant morphism is representable.
If $\XX\to Z$ is a morphism to an algebraic space $Z$
and $W\to Z$ is a morphism of algebraic spaces,
then the projection $\XX\times_ZW\to \XX$ is isovariant.
For any Artin stack $\XX$ we will denote by $\XX_{iso.et}$ the site
whose objects are isovariant \'etale morphisms $\YY\to \XX$.
A morphism from $y:\YY\to \XX$ to $y':\YY'\to \XX$ in $\XX_{iso.et}$
is a pair $(f:\YY\to \YY',\sigma)$ where $\sigma:y\simeq y'\circ f$.
Then using Quillen's Q-construction and the loop functor,
we have two presheaves of spectra
\[
\mathbf{G}:(\XX_{iso.et})^{op}\longrightarrow \mathbf{Spt}
\]
and 
\[
\mathbf{K}:(\XX_{iso.et})^{op}\longrightarrow \mathbf{Spt}
\]
which to any $\YY\to \XX$ in $\XX_{iso.et}$
associate $G(\YY)$ and $K(\YY)$ respectively,
where $\mathbf{Spt}$ is the category of spectra.
For a presheaf of spectra $P$
let $\mathbb{H}_{iso.et}(\XX,P)$ be the hypercohomology
with respect to isovariant \'etale topology on $\XX$, defined
in \cite[4.2.3]{Jo}. In loc. cit., to define hypercohomology
of presheaves of spectra the author uses Godement resolutions.
If you are familiar with model categories, you may consider
the hypercohomology of $P$ to be $P'(\XX)$ where $P\to P'$
is a fibrant replacement in the category of presheaves of
spectra endowed with the Jardine's model structure \cite[3.3]{Mi}
with respect to isovariant \'etale topology.
Similarly, for an algebraic space $X$ and any presheaf
of spectra $P$ on the \'etale site on $X$, we write
$\mathbb{H}_{et}(X,P)$ for the hypercohomology of $P$ with respect to
\'etale topology.
The descent theorem \cite[5.10]{Jo} says that
there is a weak equivalence $G(\XX)\otimes\QQ\to \mathbb{H}_{iso.et}(\XX,\mathbf{G}\otimes\QQ)$. (See \cite{Jo} for various localized versions.)
The following generalizes Poincar\'e duality \cite[5.16]{Jo},
which was proved in the case of Deligne-Mumford stacks.

\begin{Proposition}[Poincar\'e duality]
\label{duality}
Let $\XX$ be a regular Artin stack of finite type over $\ZZ$ with finite inertia stack.
Let $\pi:\XX\to X$ be a coarse moduli map and
$\pi_\#$ the direct image functor of presheaves of spectra.
Then the natural map
\[
G(\XX)\otimes\QQ\simeq \mathbb{H}_{et}(X,\pi_\#\mathbf{G}\otimes\QQ)\leftarrow \mathbb{H}_{et}(X,\pi_\#\mathbf{K}\otimes\QQ)
\]
is a weak equivalence of spectra.
\end{Proposition}

\Proof
To show our claim, clearly we may work \'etale locally on the
coarse moduli space $X$.
Thus according to Proposition~\ref{main} we may and will assume that
$\XX$ is of the form $[\Spec R/\operatorname{GL}_n]$.
Now we can apply the result of Thomason \cite[Theorem 5.7]{Theq}
to obtain our Proposition.
\QED

\begin{Remark}
Proposition~\ref{duality} also holds for other localized $G$-theories \cite[5.1.5]{Jo}.
\end{Remark}

\end{document}